\documentclass[12pt]{article}

\usepackage{amsmath}
\usepackage{multirow} 
\usepackage{latexsym,amsfonts}
\usepackage{graphicx}
\usepackage{subcaption}
\usepackage{pdflscape}
\usepackage{caption}
\usepackage{tikz-cd} 
\usepackage{tikz}
\usepackage{enumitem} 
\usepackage{authblk}
\newtheorem{thm}{Theorem}[section]

\newtheorem{prop}[thm]{Proposition}

\newtheorem{Rem}[thm]{Remark}

\author[1]{Sanja Rukavina}
\author[2]{Vladimir D. Tonchev}
\affil[1]{Faculty of Mathematics, University of Rijeka, 51000 Rijeka, Croatia}
\affil[2]{ Department of Mathematical Sciences, Michigan Technological University, Houghton, MI 49931, USA
}
\title{On symmetric $2$-$(70,24,8)$ designs with an automorphism of order $6$}

\begin{document} 
\maketitle

\begin{abstract}
In this paper we analyze possible actions of an automorphism of order six on a $2$-$(70, 24, 8)$ design, and give a complete classification for the action of the cyclic group of order six $G= \langle \rho \rangle \cong Z_6 \cong Z_2 \times Z_3$, where $\rho^3$ fixes exactly $14$ points (blocks) and $\rho^2$ fixes $4$ points (blocks). Up to isomorphism there are $3718$ such designs. This result significantly increases the number of 
previously known $2$-$(70,24,8)$ designs.
\end{abstract}


\section {Introduction}
We assume familiarity with the basic facts and notions from the theory of combinatorial designs \cite{BJL}, \cite{l}, \cite{ton88}.\\

The first symmetric $2$-$(70,24,8)$ design was constructed by Zvonimir Janko and Tran van Trung in 1984  \cite{jt70}. The Janko-Trung design is a self-dual design with a full automorphism group isomorphic to the group $Frob_{21} \times Z_2$ of order $42$. Later on, following the same method based on enumeration of cosets in a group (see \cite{jan}), A. Golemac proved that, up to isomorphism and duality, there exist five symmetric $2$-$(70, 24, 8)$ designs whose automorphism group is isomorphic to the group $E_8: Frob_{21}$ of order $168$ \cite{ag70}. The five designs that A. Golemac found are
all non-self-dual, hence there are ten nonisomorphic designs invariant under the group $E_8: Frob_{21}$ of order $168$.
In  \cite{dc70}, D. Crnkovi\'{c}  proved the existence of twenty-two nonisomorphic symmetric $2$-$(70,24,8)$ designs having a full automorphism group isomorphic to $Frob_{21} \times Z_2$. Since the designs  constructed by Crnkovi\'{c} include the first design found by Janko and Trung, the existence of $32$ symmetric $2$-$(70, 24, 8)$ designs has been established so far.
 To the best of our knowledge, there are no other successful attempts to construct $2$-$(70, 24, 8)$ designs, and no example of a $2$-$(70, 24, 8)$ design with a  full automorphism group of order smaller than 42
is known.\\

In this paper we analyze possible actions of an automorphism of order six on a $2$-$(70, 24, 8)$ design, and give a complete classification for the action of 
the cyclic automorphism group of order six $G= \langle \rho \rangle \cong Z_6 \cong Z_2 \times Z_3$ where $\rho^3$ fixes exactly $14$ points (blocks) and $\rho^2$ fixes $4$ points (blocks). For such an action, the existence of $3718$ non-isomorphic symmetric $2$-$(70, 24, 8)$ designs is proved.

\section{A construction of $2$-$(70,24,8)$ designs with an automorphism of order six}

For the construction of $2$-$(70,24,8)$ designs with an automorphism of order six we use the method for constructing orbit matrices with presumed action of an automorphism group, which are then indexed to construct designs (see, for example, \cite{c-r}, \cite{cep}, \cite{jan} ). This method is often used when a presumed automorphism group is of composite order. In particular, we use the following result.

\begin{prop} \cite[Proposition~2.3]{ord9} \label{pq}
 Let $p$ and $q$ be two distinc prime numbers and $G= \langle \rho \rangle \cong Z_{p \cdot q} \cong Z_p \times Z_q$ be a cyclic automorphism group of a symmetric block design $\mathcal{D}$, Then the $G$-orbits of points (or blocks) of the design $\mathcal{D}$ having length $p$ or $q$ consist of fixed points (or blocks) of the permutation $\rho^p$ or $\rho^q$, respectively. Furthermore, the  $G$-orbits of points (or blocks) of the design $\mathcal{D}$ having length $p\cdot q$ consist of $p$ $\langle \rho^p \rangle$-orbits of length $q$, 
 and  $q$ $\langle \rho^q \rangle$-orbits of length $p$. 
\end{prop}

Using Proposition \ref{pq}, after constructing orbit matrices for an automorphism group $G= \langle \rho \rangle \cong Z_6 \cong Z_2 \times Z_3$, we construct their refinements and obtain orbit matrices for the cyclic group 
$\langle \rho^3 \rangle \vartriangleleft G$ of order two, such that the corresponding designs admit $\rho^2$ as an automorphism. From these orbit matrices we construct symmetric designs. For a detailed explanation of the method of construction, the reader is refered to \cite{c-r}. In our work we use computers. In addition to our own computer programs, we use computer programs by V. \'{C}epuli\'{c} 
for the construction of orbit matrices and the computer algebra system MAGMA \cite{magma} when working with codes.\\

\subsection{Possible actions of an automorphism of order six on a $2$-$(70,24,8)$ design}

The first step in the construction is to determine possible orbit lengths distributions. For that we need the following results.

\begin{prop}\cite[Corollary~3.7]{l}\label{cor-fp}
Suppose that a nonidentity automorphism $\sigma$ of a symmetric $2$-$(v,k, \lambda)$ design fixes $f$ points. Then
$$f \le v-2(k-\lambda) \qquad {\rm and} \qquad f \le ( \frac{\lambda}{k- \sqrt{k-\lambda}} ) v.$$
Moreover, if equality holds in either inequality, $\sigma$ must be an involution and every non-fixed blosk contains exactly $\lambda$ fixed points.
\end{prop}

\begin{prop}\cite[Proposition~4.23]{l}\label{involution}
Suppose that ${\mathcal D}$ is a nontrivial symmetric $2$-$(v,k, \lambda)$ design with an involution $\sigma$fixing $f$ points and blocks. If $f \neq 0$ then
$$
f \ge \left \{
\begin{tabular}{l  l}
 $1 + \frac{k}{\lambda}$,    & if $k$ and $\lambda$ are both even, \\
 $1 + \frac{k-1}{\lambda}$,  & otherwise. \\
\end{tabular} \right .
$$
\end{prop}

Denote by $f_i$, $i \in \{2,3 \}$, the number of fixed points for an action of an automorphism of order $i$ on a $2$-$(70,24,8)$ design. From Proposition \ref{cor-fp} and Proposition \ref{involution} we have $f_2 \in \{0, 4, 6, 8, 10, 12, 14, 16, 18, 20, 22, 24, 26, 28 \}$ and $f_3 \in \{1, 4, 7, 10, 13, 16, 19, 22, 25\}$.\\

Suppose that an automorpishm $\rho$ of order six acts on a $2$-$(70,24,8)$ design with the orbit lengths distribution $(d_1 \times 1, d_2 \times 2,d_3 \times 3, d_6 \times 6)$, where $d_i$ denotes the number of orbits of length $i$, $i \in \{1, 2, 3, 6\}$. If $\rho^3$ fixes $f_2$ points and $\rho^2$ fixes $f_3$ points, then $d_1+2d_2=f_3$ and $d_1+3d_3=f_2$. Furthermore, $d_1+2d_2+3d_3+6d_6=70$. \\
We checked all corresponding orbit lengths distributions for an action of an automorphism group of order six $G= \langle \rho \rangle \cong Z_6 \cong Z_2 \times Z_3$, applying the method for the construction of symmetric designs from orbit matrices for presumed action of an automorphism group and Proposition \ref{pq}. The results of our analysis are given in Table \ref{OS}, where ``-'' means that corresponding orbit matrix does not exist and ``?'' means that construction of corresponding orbit matrices is out of our reach because of the large number of possibilities. If the corresponding orbit matrices are constructed for some $f_2$ and $f_3$, then ``$d$/s'' in Table \ref{OS} means that there are $d$ corresponding orbit matrices and the meaning of a string ``?''/``Y''/``N'' is ``the construction is out of our reach''/``corresponding designs exist''/``designs do not exist'', respectively.\\

\begin{table}[htpb!] 
\begin{center} \begin{footnotesize}
\begin{tabular}{|c ||c| c | c|c|c|c|c|c|c|}
 \hline   
 
$f_2 \setminus f_3$&1& 4& 7& 10&13& 16& 19& 22& 25
\cr  \hline \hline
0 & -&  1491806/$?^{*}$& - & 1871/?& - & 5/N&-&-&-\\
 \hline  
4 &  5934/?& -& -&-&-&-&-&-&- \\
 \hline 
6 & - &  141907/$?^{**}$&-&159/N&-&-&-&-&- \\
 \hline 
8 & - &  18850/?&-&15/N&-&-&-&-&- \\
 \hline 
10 & 251398/? & -&1546/N&-&-&-&-&-&- \\
\hline 
12 & - &  ?&-&239/N&-&-&-&-&- \\
\hline 
14 & - & 65205/Y&-&-&-&4/N&-&-&- \\
\hline 
16 &  ? & 739/N&87/N&-&2/N&-&-&-&- \\
\hline 
18 & - &  ?&37/N&143/N&-&-&-&-&- \\
\hline
20 & - & -&-&-&-&-&-&-&- \\
\hline
22 &  ? & 20/N &13/N&-&-&-&-&-&-\\
\hline
24 & - &  ?&-&17/N&-&-&-&-&- \\
\hline
26 & - &  ?&-&-&-&-&-&-&- \\
\hline
28 &  ? & 178/N&-&-&-&-&-&-&- \\
\hline

\end{tabular} \end{footnotesize} 
 \caption{Possible actions of an automorphism of order six on a $2$-$(70,24,8)$ design} \label{OS}
\end{center} 
\end{table}
 
The results of our analysis of possible actions of an automorphism of order six on a $2$-$(70,24,8)$ design are summarized in Theorem \ref{th_OM}.

\begin{thm} \label{th_OM}

Let $\rho$ be an automorphism of order six acting on a symmetric $2$-$(70,24,8)$ design ${\mathcal D}$. Let $f_2$ and $f_3$ denote the number of fixed points (blocks) of $\rho^3$ and $\rho^2$, respectively.

(a) If $\rho^3$ acts on ${\mathcal D}$ without fixed points (blocks), then $f_3 \in \{4, 10 \}$.

(b) If $f_2 > 0$, then $f_3 \in \{1, 4 \}$. Especially, if $\rho^2$ fixes exactly one point (block), then $f_2 \in \{4, 10, 16, 22, 28 \}$, and if $\rho^2$ fixes exactly four points (blocks), then $f_2 \in \{6, 8, 12, 14, 18, 24, 26 \}$.
 
\end{thm}

\begin{Rem}
{\rm
The construction of $2$-$(70,24,8)$ designs from orbit matrices for an action of an automorphism of order six for the cases marked with ``*'' and ``**'' in Table \ref{OS} is out of our reach. However, the existence of $2$-$(70,24,8)$ designs for these cases has been proved using an action of a group of order $42$ and $168$, respectively (for more details see \cite{dc70} for the case marked with ``*'' and \cite{ag70} for the case marked with ``**''). 
}
\end{Rem}

 \subsection{New symmetric $2$-$(70,24,8)$ designs} \label{newZ6}
 
 In our analysis described in the previous section,
we proved the existence of symmetric $2$-$(70,24,8)$ designs with the cyclic automorphism group of order six $G= \langle \rho \rangle \cong Z_6 \cong Z_2 \times Z_3$ where $\rho^3$ fixes exactly $14$ points (blocks) and $\rho^2$ fixes $4$ points (blocks).  The corresponding orbit lengths distribution is $(2 \times 1,1 \times 2,4 \times 3,9 \times 6)$ and there are $65205$ orbit matrices for that case. Further analysis of these orbit matrices shows that among them only $66$ orbit matrices produce $2$-$(70,24,8)$ designs admitting an automorphism of order six. After eliminating isomorphic copies we obtain $3718$ non-isomorphic $2$-$(70,24,8)$ designs. Table \ref{des} contains more information on the designs constructed. Note that no example of a $2$-$(70,24,8)$ design with a full automorphism group of order $6$ or $24$ has been known previously.
 All designs from Table \ref{des} are available at
\begin{verbatim}
 https://www.math.uniri.hr/~sanjar/structures/
 \end{verbatim}

\begin{table}[htpb!] 
\begin{center} \begin{footnotesize}
\begin{tabular}{|c |c| c | c|c|}
 \hline   
 
No. of & The order of & The structure of  & No. of self-dual& No. of dually \\
designs&$Aut\mathcal{(D)}$&$Aut\mathcal{(D)}3510$&designs&isomorphic pairs\\
 \hline \hline
3510& 6 & $Z_6$&10 &1750 \\
\hline
184&24&$A_4 \times Z_2$&0 &92 \\
\hline
16&42&$Frob_{21} \times Z_2$&2&7 \\
\hline
8&168&$E_8:Frob_{21}$&0&4 \\
\hline

\end{tabular} \end{footnotesize} 
 \caption{$2$-$(70,24,8)$ designs on which an automorphism of order six acts with orbit lengths distribution $(2 \times 1,1 \times 2,4 \times 3,9 \times 6)$} \label{des}
\end{center} 
\end{table}

 The results of our classification of $2$-$(70,24,8)$ designs on which an automorphism of order six acts with orbit lengths distribution $(2 \times 1,1 \times 2,4 \times 3,9 \times 6)$ are  summarized in the following theorem.
 
 \begin{thm}
 Let $\rho$ be an automorphism of order six acting on a symmetric $2$-$(70,24,8)$ design. Up to isomorphism, there are exactly $3718$ symmetric $2$-$(70,24,8)$ designs on which the group $\langle \rho \rangle$ acts so that $\rho^3$ fixes fourteen points (blocks) and $\rho^2$ fixes four points (blocks). Among these designs, there are $12$  self-dual designs and $1853$ pairs of dually nonisomorphic designs. Exactly $3510$ designs have a full automorphism group isomorphic to the cyclic group of order six. Furthermore, $184$ designs have a full automorphism group of order $24$ isomorphic to the group $A_4 \times Z_2$, $16$ designs have a full automorphism group of order $42$ isomorphic to the group $Frob_{21} \times Z_2$ and $8$ designs have a full automorphism group of order $168$ isomorphic to the group $E_8:Frob_{21}$.
 \end{thm}

 As we have already mentioned, the  classification of symmetric $2$-$(70,24,8)$ designs with a full automorphism group isomorphic to  $E_8:Frob_{21}$  was completed by A. Golemac \cite{ag70}. 
 Our designs with a full automorphism group of order $168$ coincide with the eight designs 
 with three orbits given in \cite{ag70}. The remaining two $2$-$(70,24,8)$ designs from \cite{ag70} have two orbits and they were not obtained in our construction, since in that case an involution acts with six fixed points.\\
 Furthermore, symmetric $2$-$(70,24,8)$ designs with a full automorphism group isomorphic to the group $Frob_{21} \times Z_2$ were previously constructed by D. Crnkovi\'{c} \cite{dc70}. Nine designs from \cite[Theorem 7]{dc70} did not occur in our construction, since in that case an involution acts without fixed points (blocks).  In \cite[Theorem 8]{dc70}, the existence of eight symmetric $2$-$(70,24,8)$ designs (up to isomorphism and duality) with the automorphism group $Frob_{21} \times Z_2$ acting with the orbit lengths distribution $(7,7,14,42)$ was proved. An analysis of designs given as a part of the proof of that theorem  in \cite{dc70} shows that among these designs there are two self-dual designs, which means that $14$ designs with the orbit lengths distribution $(7,7,14,42)$ were constructed in \cite{dc70}, and that the total number of designs constructed in \cite{dc70} is $23$. In our construction of designs on which an automorphism of order six acts with the orbit lengths distribution $(2 \times 1,1 \times 2,4 \times 3,9 \times 6)$, we obtained $16$ designs  with a full automorphism group isomorphic to the group $Frob_{21} \times Z_2$, and on all of them the full automorphism group acts  with the orbit lengths distribution $(7,7,14,42)$, while the subgroup isomorphic to $Frob_{21}$ acts with the  orbit lengths distribution $(7,7,7,7,21,21)$. As given in Table \ref{des}, among these designs there are two self-dual designs and seven pairs of dually isomorphic designs. Furthermore, $14$ designs are isomorphic to those given in \cite{dc70} and there is one additional pair of dually isomorphic designs not given in \cite{dc70}. Note that in the case of an action of $G \cong Frob_{21} \times Z_2$ on a $2$-$(70,24,8)$ design with the orbit lengths distribution $(7,7,14,42)$  the subgroup of $G$ isomorphic to $Z_2$ always fixes $14$ points (bloks) and $Z_3$ as a subgroup of $G$ in that case fixes $4$ points (blocks). Because of that, our construction gives a complete classification of $2$-$(70,24,8)$ designs with the automorphism group $G \cong Frob_{21} \times Z_2$ acting with the orbit lengths distribution $(7,7,14,42)$. The next theorem fixes an error in \cite[Theorem 8]{dc70}.
 
 \begin{thm}
Up to isomorphism and duality there are nine symmetric
$2$-$(70,24,8)$ designs with an automorphism group $Frob_{21} \times Z_2$ acting with the orbit lengths distribution $(7,7,14,42)$. Two of these designs are self-dual. The full automorphism groups of these designs are isomorphic to $Frob_{21} \times Z_2$.
 \end{thm}

As a consequence of our observations, we give a correction of \cite[Theorem 9]{dc70}.

\begin{thm}
 Up to isomorphism, there are $25$ symmetric
$2$-$(70,24,8)$ designs with an automorphism group isomorphic to $Frob_{21} \times Z_2$. Among them there are three self-dual designs and eleven pairs of dually isomorphic designs. The full automorphism groups of these designs are isomorphic to $Frob_{21} \times Z_2$.
\end{thm}

Our analysis shows that there are $11$ known symmetric $2$-$(70,24,8)$ designs which are not covered by our construction. Hence, the total number of known nonisomorphic designs with parameters  $2$-$(70,24,8)$
is $3729$.  Table \ref{all} 
contains information about these designs.
 
\begin{table}[htpb!] 
\begin{center} \begin{footnotesize}
\begin{tabular}{|c |c| c | c|c|}
 \hline   
 
No. of & The order of & The structure of  & No. of self-dual& No. of dually \\
designs&$Aut\mathcal{(D)}$&$Aut\mathcal{(D)}$&designs&isomorphic pairs\\
 \hline \hline
3510& 6 & $Z_6$&10 &1750 \\
\hline
184&24&$A_4 \times Z_2$&0&92 \\
\hline
25&42&$Frob_{21} \times Z_2$&3&11 \\
\hline
10&168&$E_8:Frob_{21}$&0&5 \\
\hline

\end{tabular} \end{footnotesize} 
 \caption{Known symmetric $2$-$(70,24,8)$ designs} \label{all}
\end{center} 
\end{table}

 \section {On the binary codes of $2$-$(70,24,8)$ designs }
 
 For basic facts and notions from coding theory we refer the reader to \cite{FEC}.\\
 
 In the previous section, we presented a construction of new symmetric $2$-$(70,24,8)$ designs using orbit matrices for presumed action of an automorphism of order six. However, some of the new designs can also be obtained by applying a different method, namely by analyzing the codes of known designs \cite{CRT}, \cite{MT}, which gives an interesting insight on how these designs are related.
 
 In Table \ref{2rank}, we give information about the dimensions of the binary codes spanned by the incidence matrices of designs constructed in Section \ref{newZ6}.
 The lowest 2-rank of the incidence matrix of any known 2-$(70,24,8)$ design is 22, and there are 14 designs with 2-rank 22:
 six designs with a group of order 6, six designs with a group of order 24, and two designs with a group of order 168 which
are isomorphic to  the Golemac Design $D_1$ \cite[page 57]{ag70} and its dual design $D_1^\perp$.

The binary linear code spanned by the incidence vectors of the blocks of  $D_1^\perp$ contains 49427 codewords of weight 24.
We computed the orbits of codewords of weight 24 under the action of a subgroup $H_3$ of order 3 of the automorphism group of   $D_1^\perp$
and did a complete search for 2-$(70,24,8)$ designs invariant under $H_3$. As a  result of this search exactly four distinct designs
were found:
 $D_1^\perp$ and three pairwise nonisomorphic designs $D'_1, D'_2, D'_3$, all having a full automorphism group of order 24 and
 2-rank 22.
 The designs $D'_1, D'_2, D'_3$  are not self-dual. Thus,   $D'_1, D'_2, D'_3$ and their dual designs are
  the six designs with 2-rank 22 in Table \ref{2rank}.
 
 A similar search for designs invariant under a subgroup $H_7$ of order 7 in the code of the dual design $D_2^\perp$ of  Golemac's Design $D_2$ \cite[page 57]{ag70} which is of dimension 23, shows that this code contains exactly two designs invariant under $H_7$: $D_2^\perp$ and a design isomorphic to
 Golemac's design  $D_1^\perp$.

\begin{table}[htpb!] 
\begin{center} \begin{footnotesize}
\begin{tabular}{|c ||c| c | c|c|}
 \hline   
 
 $2$-rank&No. of designs & No. of designs & No. of designs  & No. of designs  \\
 of design& $|Aut(\mathcal{D}|=6$ & $|Aut(\mathcal{D}|=24$& $|Aut(\mathcal{D}|=42$& $|Aut(\mathcal{D}|=168$ 
\cr  \hline \hline
22 & 6&  6 & 0 & 2 \\
 \hline  
23& 20 &14 & 0 &2\\
\hline
24& 308&24&2&0\\
\hline
25&604&54&2&2\\
\hline
26&1402&86&0&2\\
\hline
27&198&0&0&0\\
\hline
28&364&0&10&0\\
\hline
29&138&0&0&0\\
\hline
30&385&0&1&0\\
\hline
31&69&0&1&0\\
\hline

\end{tabular} \end{footnotesize} 
 \caption{  $2$-ranks of constructed designs} \label{2rank}
\end{center} 
\end{table}

\bigskip
\noindent {\bf Acknowledgement} 

The first author is supported by {\rm C}roatian Science Foundation under the project 6732.

\end{document}